\documentclass[a4paper,10pt]{amsart}
\usepackage[utf8]{inputenc}
\usepackage{amsmath}

\title{On Fano manifolds of Picard number one}
\author{Laurent Manivel}
\date{January 29th, 2015}

\newtheorem{theorem}{Theorem}[section]

\newtheorem{proposition}[theorem]{Proposition}

\def\PP{\mathbb{P}}

\def\ra{\rightarrow}

\begin{document}

\begin{abstract}
 K\"uchle classified the Fano fourfolds that can be obtained 
as zero loci of global sections of homogeneous vector bundles on Grassmannians. 
Surprisingly, his classification exhibits two families of fourfolds with the same discrete
invariants. Kuznetsov asked whether these two types of fourfolds are deformation equivalent. 
We show that the answer is positive in a very strong sense, since the two families are in 
fact the same! This phenomenon happens in higher dimension as well.
\end{abstract}

\maketitle

\section{Introduction}

The classification of smooth Fano threefolds by Iskhovskih, Mori and Mukai 
has been
one of the highlights of 20th century's algebraic geometry. Many interesting 
cases from that classification, and especially among the prime Fano manifolds of
index one, are obtained by taking suitable sections (mostly linear sections) of 
certain rational homogeneous spaces, and Mukai wrote a wonderful series of papers
about their astonishing geometry (see for example \cite{Mu92}, and the more general reference 
\cite{IP99}). 

It is a general fact that rational homogeneous 
spaces are a rich source of interesting Fano manifolds, obtained as zero-loci of 
global sections of 
vector bundles, and especially homogeneous vector bundles of low rank. In dimension 
four, O. K\"uchle \cite{Ku95} began the classification of these Fano manifolds by focusing 
on Fano fourfolds of index one obtained as subvarieties of Grassmannians, and 
defined as zero-loci of semisimple homogeneous vector bundles. He obtained a 
list of fourfolds which, as recently stressed by A. Iliev, are potentially a rich 
source of nice geometry. In particular some of these varieties have special 
Hodge structures, that could be relevant in the quest for new hyperk\"ahler manifolds. 
More precisely, they seem to be good candidates for the ideas of \cite{IM15} to be 
implemented successfully. 

Recently, A. Kuznetsov obtained nice structural results about the K\"uchle fourfolds
whose Picard number is bigger than one. He also observed that among those whose Picard 
group is cyclic, there are two families with the very same discrete invariants. 
He asked whether this coincidence could be explained by the possibility that the 
two types of Fano fourfolds are deformation equivalent \cite[Question 1.1]{Kuz15}. 

The main result of this short note is that this is indeed the case, and that 
much more is true: the two families are in fact the same one! Moreover, this phenomenon
happens in arbitrary dimension: there are two families of prime Fano $n$-folds of index 
one, that look different at first sight but in fact coincide. The first one is that 
of $(n+2)$-codimensional linear sections of the Grassmannian $G(2,n+3)$. The second one is that of 
zero-loci of sections of the twisted quotient bundle on $G(2,n+2)$. We prove in
Theorem \ref{main} that these two types of varieties are the same up to projective 
equivalence. Meanwhile we provide a few elements about the geometry of these Fano manifolds,
that would probably deserve further study. 

\medskip\noindent {\it Acknowledgements.} I thank A. Kuznetsov for his interesting comments. 

This work has been carried out in the framework of the Labex Archimède (ANR-11-LABX-0033) and 
of the A*MIDEX project (ANR-11-IDEX-0001-02), funded by the ``Investissements d'Avenir" 
French Government programme managed by the French National Research Agency.

\section{Two families of Fano manifolds of index one}

We will denote by $G_n$ the Grassmannian $G(2,n+2)$ parametrizing planes in a complex 
vector space $V_{n+2}$ of dimension $n+2$. This is a smooth Fano variety of
dimension $2n$, Picard number one and index $n+2$. The very ample generator of
the Picard group defines the Pl\"ucker embedding of $G_n$, with respect to which
its degree is equal to the Catalan number $c_{n}=\frac{1}{n+1}\binom{2n}{n}$.

\subsection{Linear sections of Grassmannians}

Any smooth linear section $X$ of $G_{n+1}$ of codimension $n+2$ is again a Fano variety,  
of dimension $n$, Picard number one and index one. This will be our first family
of smooth Fano manifolds. Note that for $n=2$ we get a del Pezzo surface of degree five, 
and for $n=3$ a prime Fano threefold of genus eight. 

There is a moduli space for these manifolds, that we can 
construct as the GIT quotient  of an open subset of the Grassmannian $G(n+2,\wedge^2 V_{n+3}^*)$ 
by the reductive group $PGL_{n+3}$. The dimension of this moduli space is 
$N=(n+1)(n+2)^2/2-(n+2)(n+4)=(n+3)(n^2-4)/2$. Locally around a point in the moduli space 
corresponding to a given $X$, deformations 
are unobstructed: recall this is the case for any Fano manifold, as a consequence
of the Kodaira-Akizuki-Nakano vanishing theorem. The tangent  space $H^1(TX)$ to the local 
Kuranishi space  can be computed fom the 
normal exact sequence, which yields a long exact sequence of 
cohomology groups on $X$:
$$H^0(TX)\ra H^0(TG_{n+1|X})\ra H^0(O_X(1))^{n+2}\ra H^1(TX)\ra H^1(TG_{n+1|X}).$$
Indeed, using the Koszul resolution of $O_X$ and Bott's theorem on $G_{n+1}$, one 
checks that $H^1(TG_{n+1|X})=0$ and $H^0(TG_{n+1|X})=H^0(TG_{n+1})=sl_{n+3}$. 
Moreover $H^0(TX)=0$, or equivalently:

\begin{proposition}
For $n\ge 4$, any smooth $X$ has a finite automorphism group. 
\end{proposition}

\proof Since $X$ has index one, $H^0(TX)=H^0(\Omega^{n-1}_X(1))$. Taking the $(n-1)$-th wedge 
power of the conormal exact sequence, we get that this cohomology group is zero as soon as 
$$H^k(X, \Omega^{n-1-k}_{G_{n+1}|X}(1-k))=0 \qquad\forall 0\le k\le n-1.$$
Using the Koszul resolution of the structure sheaf of $X$, we get that this vanishing will hold as soon as
$$H^{k+\ell}(G_{n+1}, \Omega^{n-1-k}_{G_{n+1}}(1-k-\ell))=0 \qquad\forall 0\le \ell\le n+2.$$
If $k+\ell\ge 2$, this follows from the Kodaira-Akizuki-Nakano vanishing theorem, since $(k+\ell )+
(n-1-k)=n-1+\ell$ is smaller than the dimension of $G_{n+1}$.  If $k+\ell\ge 1$, we get 
$H^1(\Omega^{n-1}_{G_{n+1}})$ and $H^1(\Omega^{n-2}_{G_{n+1}})$, which are both zero 
for $n\ge 4$ since   $H^q(\Omega^p_{G_{n+1}})$ is always zero for $p\ne q$. Finally if 
$k+\ell\ge 0$, we get $H^0(\Omega^{n-1}_{G_{n+1}}(1))=0$ by \cite[Theorem 2.3]{Sn86}.

\medskip\noindent {\it Remark}. It would be interesting to have closed formulas for the 
Hodge numbers of linear sections of Grassmannians. In the case we are interested in, the 
Lefschetz hyperplane theorem gives $h^{p,q}(X)=h^{p,q}(G_{n+1})$ for $p+q<n$. Hence
$h^{p,q}(X)=\delta_{p,q}\left \lceil{\frac{p+1}{2}}\right \rceil$ under this condition. 

The case $p+q=n$
is  more difficult. It would be enough to compute the holomorphic Euler characteristic
of the bundles of $p$-forms for $p\le n$, which can be done by standard techniques from
Schubert calculus but remains computationally hard. Even the topological Euler characteristic
seems difficult to compute. If we try to use the Gauss-Bonnet formula $e(X)=\int_Xc_n(TX)$
we can obtain the Chern class of the tangent bundle from the normal exact sequence. We 
deduce that if
$$P_n(x,y)=\sum_{k=0}^{n+1}p_{n,k}x^{2k}y^{2n+2-2k}=\Big[ \frac{x^{n+2}}{(1+x)^{n+2}}
   \frac{(1+x+y^2)^{n+3}}{1-x^2+4y^2} \Big]_{2n+2}$$
the degree $2n+2$ part of the Taylor expansion of this rational function, 
then 
$$e(X)=\sum_{k=0}^{n+1}p_{n,k}c_k.$$
We can deduce the missing Betti number $b_n=b_n(X)$ for $X$ of small dimension $n$:
 $b_2=5$,   $b_3=10$,   $b_4=69$,   $b_5=380$,   $b_6=2321$,   $b_7=9442$.

\subsection{Zero loci of the twisted quotient bundle}

Recall that the Grassmannian $G_n$ is endowed with two natural vector bundles, 
the tautological rank two bundle $U$, and the quotient bundle $Q$, which has rank $n$. 
Moreover $\det(U^*)=\det(Q)=O_{G_n}(1)$, the very ample generator of the Picard group. 
The quotient bundle $Q$ is generated by global sections, and the zero locus of a non 
zero section is just a projective space. More interesting are the zero loci of global
sections of the twisted quotient bundle 
$Q(1)$. Since $\det(Q(1))=O_{G_n}(n+1)$, the zero locus $Y$ of a general 
global section of $Q(1)$ is a smooth Fano manifold of dimension $n$ and index one. 
This is our second family of such manifolds. 

In classical langage, $Y$ defines a congruence of lines in $\PP (V_{n+2})=\PP^{n+1}$. 
Recall that the {\it order} of such a congruence is defined as the number of lines from $Y$ passing 
through a general point in $\PP^{n+1}$.

\begin{proposition}
The congruence of lines defined by  $Y$ has order $n+1$.
\end{proposition}

\proof The set of lines passing through a point in $\PP^{n+1}$ is isomorphic to 
$\PP^{n}$, and the restriction of $Q(1)$ to this projective space is $R(1)$, if
$R$ denotes the tautological quotient bundle on $\PP^{n}$. The order of the 
congruence is the number of zeroes of the induced section of $R(1)=T\PP^{n}$. 
For $Y$ general and a general point in $\PP^{n+1}$ we get a general vector 
field on $\PP^{n}$, and we deduce that the order is $c_n(T\PP^{n})=e(\PP^{n})=n+1$.
$\Box$

\medskip The space of global sections of $Q(1)=Q\otimes \det(U^*)$ is
$$S_n:= S_{10\cdots 0-1-1}V_{n+2}=Ker(V_{n+2}\otimes\wedge^2V_{n+2}^*\longrightarrow V_{n+2}^*)$$
(the latter morphism being the natural contraction map), as follows from the Borel-Weil
theorem. Its dimension is $n(n+2)(n+3)/2$. 

We would be tempted to think of our family of Fano manifolds $Y$,  
as the quotient of an open subset of $\PP (S_n)$ by $PGL_{n+2}$. 
This is not correct. Indeed, an unusual phenomenon happens: $H^0(TG_{n|Y})$
is bigger than the $sl_{n+2}$ we would have expected. In fact, 
$$H^0(TG_{n|Y})=sl_{n+2}\oplus V_{n+2}.$$
This means that there are more linear isomorphisms beteween these varieties than those coming
from $PSL_{n+2}$. We will explain the appearance of that extra factor in the next section.

\medskip
Note that $Q(1)=Hom(\wedge^2U, Q)$, and  $S_n\subset Hom(\wedge^2V_{n+2},V_{n+2})$.
Hence, for any $\omega\in S_n$, the zero locus of the associated 
section of $Q(1)$ is 
$$Y_\omega=\{\langle a,b\rangle\in G(2,V_{n+2}), \quad \omega(a,b)\in \langle a,b\rangle \}.$$
Note that this makes sense for any $\omega \in Hom(\wedge^2V_{n+2},V_{n+2})=S_n
\oplus V_{n+2}^*$; but the component on $V_{n+2}^*$ is in fact unsignificant, since for any  $v^*\in V_{n+2}^*$,
$v^*(a,b)=v^*(b)a-v^*(a)b$ always belongs to $ \langle a,b\rangle $.

More interestingly, the previous description shows that $Y_\omega$ has a natural rational map 
to $\PP (V_{n+2})$, that we denote by $\Omega$. By definition 
$$\Omega(\langle a,b\rangle)=[\omega(a,b)].$$

\begin{proposition}
The rational map $\Omega$ is a birational isomorphism with a determinantal hypersurface of 
degree $n+1$ in $\PP^{n+1}$. 
\end{proposition}

\proof Let $Z=\Omega (Y)$. A point $[c]$ of $\PP (V_{n+2})$ belongs to $Z$ if and only if 
there exists an independant vector $d$ such that $\omega(c,d)=0$ belongs to $c^{\perp}$. 
Otherwise said,  the induced map from $V_{n+2}/\langle c\rangle\ra (c^{\perp})^*$ must not be 
injective. Note that there is a natural duality between $c^{\perp}$ and $V_{n+2}/\langle c\rangle$. 
Globally, we can therefore describe $Z$ as the first
degeneracy locus of a morphism between vector bundles
$$\bar{\omega}: R(-1)\longrightarrow R.$$
This implies that $[Z]=c_1(Hom(R(-1),R))$ is $n+1$ times the hyperplane class, so $Z$ is 
a hypersurface of degree $n+1$. Moreover the 
fiber of $\sigma$ over any point can be identified with the kernel of $\bar{\omega}$, in 
particular it is always a linear space, and it reduces to a single point for the general 
point of $Z$. 

Let $\bar{Y}\subset F(1,2,V_{n+2})$ be the variety parametrizing  incident lines and planes $ 
\langle c\rangle\subset \langle c,d\rangle$ such that $\omega(c,d)$ belongs to $\langle c\rangle$. 
The two projections yield  a diagram

$$\begin{array}{rcccl}
& & \bar{Y} & & \\
& { }^{p_1}\!\textstyle{\swarrow} & &  \searrow \!\! { }^{p_2} & \\
 G(2,V_{n+2})\supset Y & & &  & Z\subset \PP (V_{n+2}) 
\end{array}$$

%\xymatrix{
% & \bar{Y} \ar@{->}[rd] \ar@{->}[ld] & \\
% G(2,V_{n+2})\supset Y & \ar@{->}[rr]  & Z\subset \PP (V_{n+2}) }

It is easy to see that for $\omega$ generic, the variety $\bar{Y} $ is the zero-locus of a generic section 
of a globally generated vector bundle, hence a smooth variety. Moreover its projection to $Y$ is birational. 
More precisely, this projection is the blow-up of the smooth  subvariety $S_\omega$ defined as 
$$S_\omega=\{\langle a,b\rangle\in G(2,V_{n+2}), \quad \omega(a,b)=0 \}.$$
Note that $S_\omega$ is a Calabi-Yau variety, of codimension two in $Y$.

The projection to $Z$ fails to be an isomorphism over the locus where $\bar{\omega}$ drops rank. 
For $\omega$ generic, this occurs on a codimension three subvariety $C\subset Z$ which is also the singular locus of 
$Z$. The projection from $\bar{Y} $ to $Z$ has fibers over the general points of $C$ which are lines, in particular 
this projection is a small morphism. $\Box$

\medskip\noindent {\it Remark}. For $n=2$, the surface $Y$ is a del Pezzo surface of degree five
and $\bar{Y}$ is its blow-up at two points. The projection to the cubic surface $Z$ is an isomorphism. 

For $n=3$, the threefold $Y$ is a prime Fano threefold of genus eight, 
and $\bar{Y}$ is obtained by blowing-up the elliptic curve $E=S_\omega$. The projection to the quartic 
determinantal threefold $Z$ contracts $25$ lines to the $25$ singular points of $Z$. 

\medskip\noindent {\it Remark}. As Kuznetsov points it out,  $\omega\in Hom(\wedge^2V_{n+2},
V_{n+2})$ might be considered as defining a bracket $[a,b]=\omega(a,b)$, although not a 
Lie bracket in general since the Jacobi identity has no reason to hold. Then $Y_\omega$
parametrizes the planes in $V_{n+2}$ on restriction to which the bracket defines a Lie 
algebra structure: each plane in $Y_\omega$ is required to be stable under the bracket, and the Jacobi identity
automatically holds for dimensional reasons. Moreover the codimension two subvariety $S_\omega$
can be interpreted as parametrizing the two-dimensional abelian subalgebras. 

\section{And their coincidence}

Let $X$ be a smooth linear section of $G_{n+1}$, defined by an $(n+2)$-dimensional space
of linear forms $H_{n+2}\subset \wedge^2V_{n+3}^*$. We suppose in this section that $H_{n+2}$
is generic. 

Fix a decomposition $V_{n+3}=V_{n+2}
\oplus \langle v_{n+3}\rangle$, yielding a decomposition $$\wedge^2V_{n+3}^*=\wedge^2V_{n+2}^*\oplus 
v_{n+3}^*\wedge V_{n+2}^*,$$ 
where the linear form $v_{n+3}^*$ has kernel $V_{n+2}$. 
Generically, the projection on the second factor 
of this decomposition yields an
isomorphism $H_{n+2}\simeq V_{n+2}^*$. The variety $X$ is thus defined by a
monomorphism $\omega\in Hom(V_{n+2}^*, \wedge^2V_{n+2}^*)$ (we use the same notation for $\omega$
and its transpose): it is cut out by the space of 
linear forms defined as the graph
$$H_{n+2}=\{\omega(u)+v_{n+3}^*\wedge u, \quad u\in V_{n+2}^*\}.$$
Since $\omega$ is injective, $X$ does not contain any line passing through $[v_{n+3}]$. 
Any line in $X$ is of the form $\langle a,b+\chi(b)v_{n+3}\rangle$ for some non zero vectors
$a,b\in V_{n+2}$, with the condition that
$$\chi(b)u(a)=\omega(u)(a,b) \quad \forall  u\in V_{n+2}^*.$$
This implies that $\omega(u)(a,b)=0$ for all  $u\in a^\perp$, while the remaining 
equation determines $\chi(b)$. This means in particular that the linear projection 
from $\wedge^2V_{n+3}$ to $\wedge^2V_{n+2}$, which  induces a rational map 
$$G_{n+1}=G(2,V_{n+3}) \dashrightarrow G_n=G(2,V_{n+2})\simeq G(2,V_{n+3}/\langle v_{n+3}\rangle),$$
restricts to a well-defined map from $X$ 
to the subvariety $Y_\omega$ of $G_n$, which is moreover injective. 

The inverse mapping can be described as follows. Note that $\omega(a,b)$ is
the vector $c\in V_{n+2}$ defined by the identity $u(c)=\omega(u)(a,b)$ for all $u\in  V_{n+2}^*$.
The equations defining $X$ reduce to the condition that $c=\chi(b)a$. 
The line $\langle a,b+\chi(b)v_{n+3}\rangle$ is thus represented by 
$$[a\wedge (b+\chi(b)v_{n+3})]=[a\wedge b+c\wedge v_{n+3}]=[a\wedge b+\omega(a,b)\wedge v_{n+3}].$$
This concludes the proof of our main result:

\begin{theorem}\label{main}
$X\subset G_{n+1}$ and $Y_\omega\subset G_n$ are projectively equivalent.
\end{theorem}

%Note that we can change the decomposition of $V_{n+3}$ to the effect of changing 
%$v_{n+3}^*$ into $v_{n+3}^*-v^*$ for some $v^*\in V_{n+2}^*$. Then $\omega(u)$ becomes
%$\omega(u)+v^*
\medskip\noindent {\it Remark}. In order to identify $X$ with the subvariety $Y_\omega$ of $G_n$, 
we started from a decomposition $V_{n+3}=V_{n+2}\oplus \langle v_{n+3}\rangle$. Note that if we choose another
hyperplane $V_{n+2}$, or equivalently if we change the defining linear form  $v_{n+3}^*$ into 
$v_{n+3}^*-e^*$ for some $e^*\in V_{n+2}^*$, then $\omega$ is changed into the morphism from 
$V_{n+2}^*$ to $\wedge^2V_{n+2}^*$ that sends $u$ to $\omega(u)+e^*\wedge u$. In particular 
the class of $\omega$ in $S_n$ is not affected. 

On the contrary, changing the line $\langle v_{n+3}\rangle$ has a non trivial effect on 
the class of $\omega$ (that could easily be expressed explicitely), but does not affect
$Y=Y_\omega$ up to projective equivalence. This is precisely what explains the extra factor 
$V_{n+2}$ inside $H^0(TG_{n|Y})$.

\medskip\noindent {\it Remark}. The Calabi-Yau two codimensional subvariety $S$ of $Y$ can
be seen directly in $X$ as the intersection of $G(2,V_{n+2})\subset G(2,V_{n+3})$ with the linear space 
that defines $X$. Of course there is a whole family of such Calabi-Yau`s in $Y$, parametrized by
an open subset of the projective space of hyperplanes in $V_{n+3}$. In particular, for 
$n=2$ we get a four dimensional family of K3 surfaces covering $Y$. 

\medskip\noindent {\it Question}. As we have seen, our two families of Fano manifolds of index 
one coincide generically. An intriguing question is to decide whether they coincide stricto sensu:
can any smooth member of each family be described as a member of the second family? A negative 
answer would be particularly interesting, as a new example of the pathological behaviors of
the moduli spaces of Fano varieties.

\bigskip 

{\scriptsize
{\sc Institut de Mathématiques de Marseille,  UMR 7373 CNRS/Aix-Marseille Universit\'e, 
Technop\^ole Château-Gombert, 
39 rue Frédéric Joliot-Curie,
13453 MARSEILLE Cedex 13,
France}

{\it Email address}:  {\tt laurent.manivel@math.cnrs.fr}
}

\end{document}